\documentclass[a4paper,11pt, reqno]{amsart}
\usepackage{latexsym}
\usepackage{amsbsy}
\usepackage{amsfonts}
\usepackage{amsmath}
\usepackage{amssymb}
\usepackage{geometry}
\usepackage[T1]{fontenc}

\usepackage{graphics,graphicx}
\usepackage [dvips]{epsfig}

\def\marginpar#1{\ignorespaces}

\newtheorem{thm}{Theorem}

\theoremstyle{definition}

\newtheorem{example}[thm]{Example}

\numberwithin{equation}{section} 
\numberwithin{thm}{section}

\makeatother
\begin{document}

\date{\today}

\title[Propagation of Epistemic Uncertainty]{Propagation of Epistemic Uncertainty in Queueing Models with Unreliable Server using Chaos expansions}

\author[K. BACHI]{{Katia} Bachi}
\address{{\bf {Katia} BACHI}\\  Operational Research Department, Research Unit LaMOS (Modeling and Optimization of Systems), 
Faculty of Exact Sciences, Campus of Targua Ouzemour, Bejaia 06000, Algeria.}
\email{bachi.katia13@gmail.com}

\author[C. CHAUVIERE]{ {C\'edric} Chauvi\`ere}
\address{{\bf {C\'edric} CHAUVI\`ERE}\\ Laboratoire de Math\'ematiques Blaise Pascal (LMBP), CNRS UMR 6620, Universit\'e Clermont Auvergne, 
Campus Universitaire des C\'ezeaux, 3 place Vasarely, TSA 60026, CS 60026, 63 178 Aubi\`ere Cedex France.}
\email{Cedric.Chauviere@uca.fr}

\author[H. Djellout]{ {Hac\`ene} Djellout }
\address{{\bf {Hac\`ene} DJELLOUT}\\ Laboratoire de Math\'ematiques Blaise Pascal (LMBP), CNRS UMR 6620,
Universit\'e Clermont Auvergne, Campus Universitaire des C\'ezeaux, 3 place Vasarely,
TSA 60026, CS 60026, 63 178 Aubi\`ere Cedex France.}
\email{Hacene.Djellout@uca.fr}

\author[K. ABBAS]{ {Karim} ABBAS}
\address{{\bf {Karim} ABBAS}\\ Operational Research Department, Research Unit LaMOS (Modeling and Optimization of Systems), 
Faculty of Exact Sciences, Campus of Targua Ouzemour, Bejaia 06000, Algeria.}
\email{kabbas.dz@gmail.com}

\begin{abstract}
In this paper, we develop a numerical approach based on Chaos expansions to
analyze the sensitivity and the propagation of epistemic uncertainty through
a queueing systems with breakdowns. Here, the quantity of interest is the
stationary distribution of the model, which is a function of uncertain
parameters. Polynomial chaos provide an efficient alternative to more
traditional Monte Carlo simulations for modelling the propagation of
uncertainty arising from those parameters. Furthermore, Polynomial chaos
expansion affords a natural framework for computing Sobol' indices. Such
indices give reliable information on the relative importance of each
uncertain entry parameters. Numerical results show the benefit of using Polynomial
Chaos over standard Monte-Carlo simulations, when considering statistical moments 
and Sobol' indices as output quantities.

\end{abstract}
\maketitle

\vspace{2pt}

\section{Introduction}

Queueing systems with breakdowns are widely used to model problems occuring in computer, manufacturing systems, 
and communication networks. Typically, in queueing models the performances measurements are
assessed for deterministic parameters values. However, in practice, the
exact values of these parameters are not well known (uncertain), and they
are generally estimated empirically from few experimental observations.
This lack of precise information propagates uncertainties in the output
measure and such study is the purpose of this work.
\medskip

Since the pioneering work of Thiruvengadam \cite{Thiruvengadam} and
Avi-Itzhak and Naor \cite{Avi}, many authors have been investigating the
queueing systems with server breakdowns, see for example \cite
{Cao82,Li97,Wang01} and references therein. In this paper, we choose a
functional approach to the analysis of the dependence of performance
measures of certain queues, such as the M/M/1 and the M/G/1 queue with
breakdowns, with respect to some input parameters. More specifically,
denoting the probability of a server breakdown by $\theta $, we seek to
compute the stationary distribution of the queue-length process, denoted by $
\pi _{\theta }$. For a fixed value of $\theta $ and a finite waiting
capacity, $\pi _{\theta }$ can be found numerically by solving $\pi _{\theta
}\Xi_{\theta }=\pi _{\theta }$ and $\sum \pi _{\theta }=1$, where $\Xi_{\theta }$
denotes the transition probability matrix of an embedded jump chain. The
definition of the embedded jump chain will depend on the type of the queue:
for the M/M/1 queue we use the sample-chain, embedded at appropriate
Poisson-times, and for the M/G/1 queue we will embed the chain at departure
and repair moments; details will be provided later in the paper. In case of a
large or infinite waiting capacity, $\pi _{\theta }$ can be obtained via a
Laplace-Stieltjes transform \cite{abate}. Solving for $
\pi _{\theta }$ involves numerical inversion, which can be computationally
demanding \cite{shortle}. The computation of $\pi _{\theta }$ is a challenging
problem and a variety of approaches have been proposed in the literature for
approximately or indirectly solving the stationary distribution. The
predominant approach is to obtain either the generating function of $\pi
_{\theta }$ or an analytical expression for $\pi _{\theta }$ containing a
Laplace-Stieltjes transform, see, for example, \cite{abate,Baccelli}. Also
numerical solutions by means of the matrix geometric method \cite{Neuts} are
available, see \cite{Mitrany,Neuts79} for details.

\medskip

In performance analysis, one is not only interested in evaluating the system
for specific set of parameters but also in the sensitivity of the performance
with respect to these parameters. For example, in a queueing model with
breakdowns, the breakdown probability is a parameter of key interest and, in
this paper, we will analyze the dependence of $\pi _{\theta }$ on $\theta $,
which is significantly more challenging than evaluating $\pi _{\theta }$ for
a fixed $\theta $. Most of the time, it is assumed that these stochastic
models are solved for fixed parameters values. However, the parameters of
the model are determined through insufficient statistical data (a limited
number of observations), leading to uncertainty in the assessment of their
values. This parametric uncertainty, induced from the incomplete information
of the parameter, is called epistemic uncertainty \cite
{Philip,Mishra1,Robert}.

\medskip

In order to estimate the uncertainty of the parameters in the performance
measures of the model, two complementary approaches may be used: uncertainty
analysis and sensitivity analysis. Uncertainty analysis consists in
modelling input parameters as random variables, see, for example, \cite
{Cacuci,Helton94,Helton97}. Then, the sensitivity analysis aims at
determining the relative contribution of individual parameters. Several
approaches for the propagation of uncertainty have been developed, including
interval arithmetic \cite{Ramon,Rocco}, Taylor series expansion \cite
{Dhople2012,Granger,sofiane2015,Shooman,Takhedmit}, moments \cite
{Granger,Padulo,Shooman}, Monte Carlo analysis \cite{Granger,Padulo,Shooman}
. Overviews of these approaches are available in several reviews \cite
{Cacuci04,Helton93,Helton06,Iman,Ionescu,Ronen}.
\medskip

Since the seminal work of Ghanem an Spanos in the 90's \cite{Gha1991} and
later generalized by Xiu and Karniadakis \cite{Xiu2002}, Polynomial Chaos
has spread over a broad scientific community, see \cite{Sepahvand}, \cite{chauviere}. 
Then, due to the specific form of these polynomials, it quickly
became apparent that this representation was suitable to compute Sobol'
indices \cite{sudret}. To the best of our knowledge, this is the first
atempt to use chaos expansion to investigate the propagation of the
epistemic uncertainty through the performance measure of the queueing
systems.
\medskip

The main objective of this work is to develop a numerical procedure to
investigate the sensitivity and the propagation of epistemic uncertainty in
the input parameters of a queueing systems with unreliable server. As a
typical example, a customer arrives at the queueing system according to a
Poisson process with rate $\lambda $, which we consider as an input
parameter of the model. We assume that this parameter is not precisely known
and therefore it can be model as a random variable of known probability
density function. Clearly, the output measure is a function of this input
random variable. The epistemic uncertainty is propagated through a
functional relationship of the type $Y(\omega )=g(\lambda (\omega ))$ that links
the input parameter to the output measures. In this work, special attention
is given to the output stationary distribution. From a statistical
point of view the quantity of interest may be the moments of the functional stationary distribution. 
Such quantities can easily be
computed from the coefficient of the polynomial Chaos representation of $
g(\lambda (\omega ))$.

Over the past few years, several techniques have emerged to compute the
coefficients of polynomial Chaoses. They can be classified in two main
families: intrusives methods  and non-intrusive methods \cite{Pettersson}. In our
context, $Y(\omega )$ can be explicitely expressed as a functional
relationship and therefore the two approaches are of equal complexity. In
this paper, projection method based on the Gauss quadrature rules are used
to obtain the chaos coefficients. Then statistical moments and Sobol'
sensitivity indices can easily be computed from polynomials chaos expansion.

\medskip

The paper is organized as follows. The next section, is devoted to restate
some basic notions about orthogonal polynomial and Chaos expansion. Then,
the Sobol' indices which are an important tool in sensitivity analysis are
introduced. Epistemic uncertainty in input parameters of two queueing
systems with unreliable server is considered in Section \ref{sec4}. In the
first queue, a single input parameter is considered as random. On the other
hand, the second queueing system considers four input random variables and
a sensitivity analysis is performed. For the two queues, numerical results
for the computation of statistical moments are obtained using both chaos
expansion and Monte-Carlo simulation. Finally, Section \ref{sec5} draws some
conclusions.


\section{Polynomial Chaos Expansion}

\label{sec1} 

Polynomial Chaos expansion shares many elements with power series involving
orthogonal polynomials, which are commonly used in the spectral community.
This is why the following subsection is devoted to restate some results
concerning orthogonal polynomials.

\subsection{Orthogonal Polynomials: univariate case}

$\quad$

Let $V$ be the real vector space of all polynomials in one variable with
real coefficients and let be the positive inner product on $V$ defined as 
\begin{equation*}
<u,v>=\int_{I}u(x)v(x)f(x)dx\quad\text{ }\forall u,v\in V
\end{equation*}
where $f:I\subset \mathbb{R}\rightarrow \mathbb{R}^{+}$ is a nonnegative
integrable function of $x$ hereinafter referred to as a `weight function'.
The set of polynomials $\{\Psi _{n}\}_{n\geq 0}$ are said to be orthogonal
with respect to the weight function $f(x)$, if $\Psi_{n}$ is a polynomials
of degree $n$ and 
\begin{equation}
<\Psi _{n},\Psi _{m}>=\int_{I}\Psi _{n}(x)\Psi _{m}(x)f(x)dx=h_{n}^{2}\delta
_{n,m}~,\qquad n,m\in \mathbb{N}\mathbb{,}  \label{eq20}
\end{equation}
where $\delta _{n,m}$ is the Kronecker's delta function, and $h_{n}$ are
non-zero constants. We recall that, for orthogonal polynomials of degree $
d=0 $, $\Psi _{0}$ is always equal to one ($\Psi _{0}=1)$. Furthermore, the
system (\ref{eq20}) is called orthonormal if $h_{n}=1$.

The most general way to build such polynomials rely on the three-term
recurrence relation

\begin{equation}
\left\{ 
\begin{array}{l}
\Psi _{n+1}(x)=(x-a_{n})\Psi _{n}(x)-b_{n}\Psi _{n-1}(x) \\ 
\Psi _{0}(x)=1,\text{ \ }\Psi _{-1}(x)=0
\end{array}
\right.  \label{eq40}
\end{equation}
with 
\begin{equation}
\left\{ 
\begin{array}{l}
a_{n}=\dfrac{<x\Psi _{n},\Psi _{n}>}{<\Psi _{n},\Psi _{n}>},\text{ }n\in 
\mathbb{N} \\ 
\\
b_{n}=\dfrac{<\Psi _{n},\Psi _{n}>}{<\Psi _{n-1},\Psi _{n-1}>},\text{ }n\in 
\mathbb{N}^{\ast }
\end{array}
\right. .  \label{eq50}
\end{equation}
(see Gautschi \cite{Gautschi} for more details). Furthermore, according to
Favard's theorem, for a given weight function $f$, corresponds a unique set
of coefficients $(a_{n},b_{n})_{n\in \mathbb{N}}$. From a numerical point of
view, the integrals appearing into (\ref{eq50}) can be evaluated using Fej\'{e}r quadrature rule, see \cite{Rahman}, for example.
Then, taking an orthogonal polynomial basis of degree $n$, $\{\Psi _{0},\Psi
_{1},\ldots,\Psi _{n}\},$ any sufficiently regular function $u:I\subset \mathbb{
\mathbb{R}}\rightarrow \mathbb{R}$ may be represented by its projection $
\Pi_{n}u$ on such basis, i.e.
\begin{equation}
\Pi_{n}u(x)=\sum_{i=0}^{n}\widehat{u}_{i}\Psi _{i}(x),  \label{eq53}
\end{equation}
where the coefficients of the projection can be computed by evaluating $
\widehat{u}_{i}=\frac{<u,\Psi _{i}>}{<\Psi _{i},\Psi _{i}>}$. Depending on
the regularity of the function $u$ and the choice of the polynomial basis,
upper bounds of the truncation error $\left\Vert u-\Pi_{n}u\right\Vert $ can
be derived \cite{Gottlieb}. In particular, when $u$ is $C^{\infty }$, one
can observe the so-called spectral convergence (the truncation error decays
exponentially fast with respect to $n$).

The construction of multivariate orthogonal polynomial basis rely on
univariate polynomial basis as will be described in the following section.

\subsection{Orthogonal Polynomials: multivariate case}

$\quad$

The setting is now the real vector space $V^{d}$ of all polynomials in $d$
variables with real coefficients, together with the positive inner product
on $V^{d}$ defined as

\begin{equation}
<u,v>=\int_{I}u(\mathbf{x})v(\mathbf{x})f(\mathbf{x})d\mathbf{x}\text{ \ }
\forall u,v\in V^{d},  \label{eq60}
\end{equation}
where, now, $f:I\subset \mathbb{\mathbb{R}}^{d}\rightarrow \mathbb{\mathbb{R}
}^{+}$ is a nonnegative integrable function of $\mathbf{x=}
(x_{1,}x_{2,}\ldots,x_{d})$. We further assume that the weight function $f$ can
be written as the product of univariate weight functions, i.e. $
f(x_{1,}x_{2,}..,x_{d})=f_{1}(x_{1})f_{2}(x_{2})\ldots f_{d}(x_{d})$, where $
f_{i}:I_{i}\subset \mathbb{\mathbb{R}}\rightarrow \mathbb{\mathbb{R}}^{+}$
are nonnegative integrable function. It would be tempting to construct the
elements of the multivariate polynomial basis in a similar fashion i.e. as a
tensor product of monovariate elements. However, such construction leads to
uncomplete basis (all monomials up to a certain degrees are not
represented). To ensure that $\{\Psi _{0}(\mathbf{x}),\Psi _{1}(\mathbf{x}
),..,\Psi _{P_{d}^{p}}(\mathbf{x})\}$ is a multivariate polynomial basis of
degree $d$ with $(P_{d}^{p}+1)$ elements, the procedure described in the sequel is usually adopted. We
first define $\Psi _{\alpha _{1}\ldots\alpha _{d}}(\mathbf{x})$ as the tensor
product of the elements of a univariate polynomial basis, of degree $n$ i.e. 
\begin{equation}
\Psi _{\alpha _{1}\ldots\alpha _{d}}(\mathbf{x})=\prod_{i=1}^{d}
\Psi_{\alpha_{i}}(x_{i}),  \label{eq70}
\end{equation}
where $\alpha _{i}\in \{0,1,\ldots,n\}$. However, not all the elements of the
form (\ref{eq70}) are retained when constructing the multivariate polynomial
basis (that would lead to an uncomplete basis of degree $nd$ with $(n+1)^{d}$
elements). Instead, for a given degree $p$, only the elements that satisfy $
\displaystyle
\sum_{i=0}^{d}\alpha _{i}\leq p$ in (\ref{eq70}) are kept and a one to one
correspondence between the multi-index $(\alpha _{0},\ldots,\alpha _{d})$ and
the $i^{th}$ element $\Psi _{i}(\mathbf{x})$ of the multivariate basis is
set. Proceeding that way, it can be shown that the number of elements of a
complete multivariate polynomial basis of degree $p$ is
\begin{equation}
P_{d}^{p}+1=\left( 
\begin{array}{c}
p+d \\ 
d
\end{array}
\right) =\frac{(p+d)!}{d!p!}.  \label{eq80}
\end{equation}

Similarly to (\ref{eq53}), any function $u:\mathbb{\mathbb{R}}
^{d}\rightarrow \mathbb{\mathbb{R}}$ may be represented by its projection $
\Pi_{p}u$ on such basis, i.e.
\begin{equation}
\Pi_{p}u(\mathbf{x})=\sum_{i=0}^{P_{d}^{p}}\widehat{u}_{i}\Psi _{i}(\mathbf{x}
). \label{eq83}
\end{equation}

Having set the basic framework of orthogonal polynomials, we now look into
the way they can be efficiently used in the field of probability, where they
are usualy refered as Polynomial Chaos (PC).

\subsection{Polynomial Chaos expansion}

$\quad$

Denote by $(\Omega ,\mathcal{A},\mathcal{P})$ the probability space, where
as usual $\Omega $ is the set of all possible outcomes, $\mathcal{A}$ is a $
\sigma $-algebra over $\Omega ,$ and $\mathcal{P}$ is a function ${\mathcal{\
A}\rightarrow [0,1]}$ that gives a probability measure on $\mathcal{A}$.
Consider an $\mathbb{R}^{d}$-valued independent random vector $\mathbf{X}
=(X_{1},\cdots,X_{d})$ that describes input uncertainties. The probability law
of $\mathbf{X}$ may be defined by the probability density function 
\begin{equation}
f_{\mathbf{X}}(\mathbf{x})=\prod\limits_{i=1}^{d}f_{i}(x_{i}),  \label{eq90}
\end{equation}
where $f_{i}(x_{i})$ is the marginal probability density of $X_{i}$ defined
on $(\Omega _{i},\mathcal{A}_{i},\mathcal{P}_{i})$. Let us now denote $
L_{2}(\Omega _{i},\mathcal{A}_{i},\mathcal{P}_{i})$ the space of real random
variables with finite second order moments, i.e. such that 
\begin{equation}
E(X_{i}^{2})=\int x_{i}^{2}f_{i}(x_{i})dx_{i}<\infty ,  \label{eq100}
\end{equation}
where $E$ stands for the mathematical expectation. $L_{2}(\Omega _{i}, 
\mathcal{A}_{i},\mathcal{P}_{i})$ is a Hilbert space which can be provided
with a set of complete orthogonal basis $\{\Psi _{j}^{i}(x)\}_{j\ge 0}$
that are consistent with the density of $X_{i}$. For example, Legendre
polynomials are associated with uniform distributions whereas Hermite
polynomials are associated with Gaussian distributions. Similarly, $
L_{2}(\Omega ,\mathcal{A},\mathcal{P})$ is provided with a set of complete
multivariate orthogonal basis $\{\Psi _{j}(\mathbf{x})\}_{j\ge 0}$ which,
in turn, is consistent with the density of $\mathbf{X}$.

Let $\mathbf{Y}=(Y_{1},\ldots,Y_{\ell}):\Omega \rightarrow \mathbb{R}^{\ell}$ such
that $Y_{i}\in L_{2}(\Omega ,\mathcal{A},\mathcal{P})$ for $i=1,\ldots,\ell$ be a
mathematical model. For the sake of simplicity, we will consider only one
component of this model, denoted by $Y$, since the same procedure apply
identically to all the other components. Since $Y$ is assumed to belong to $
L_{2}(\Omega ,\mathcal{A},\mathcal{P})$, it can be represented as \cite{Gha1991,Xiu2002}
\begin{eqnarray}
Y(X_{1},..,X_{d}) &=&z_{0}\Psi _{0}+\sum\limits_{i_{1}=1}^{\infty
}z_{i_{1}}\Psi _{1}(X_{i_{1}})+\sum\limits_{i_{1}=1}^{\infty
}\sum\limits_{i_{2}=1}^{i_{1}}z_{i_{1}i_{2}}\Psi _{2}(X_{i_{1}},X_{i_{2}})+ 
\notag \\
&&+\sum\limits_{i_{1}=1}^{\infty
}\sum\limits_{i_{2}=1}^{i_{1}}\sum
\limits_{i_{3}=1}^{i_{2}}z_{i_{1}i_{2}i_{3}}\Psi
_{3}(X_{i_{1}},X_{i_{2}},X_{i_{3}})+\cdots,  \label{eq110}
\end{eqnarray}
which, after some rearranging, can be rewritten in a more convenient way as 
\begin{equation}
Y(\mathbf{X})=\sum\limits_{j=0}^{\infty }y_{j}\Psi _{j}(\mathbf{X})
\label{eq120}
\end{equation}
Then, similarly to equation (\ref{eq83}), this serie is truncated by keeping
terms up to a degree $p$ 
\begin{equation}
Y(\mathbf{X})\approx \Pi_{p}Y(\mathbf{X})=\sum\limits_{j=0}^{P_{d}^{p}}y_{j}
\Psi _{j}(\mathbf{X}).  \label{eq130}
\end{equation}

\begin{example}\label{ex1}$\quad$

Let $X=(X_{i_{1}},X_{i_{2}})$ be a random normal vector.
The functional approximation of the response $Y=Y(\mathbf{X})$ may be
approximated with Hermite-Chaos expansions according to (\ref{eq130}). If we
fix the degree of the Chaos expansions to be $d=2$, we obtain the series
expansion with $P_{2}^{2}=6$ terms as follows: 
\begin{eqnarray}
Y &\approx &\sum_{i=0}^{5}y_{i}\Psi _{i}(X_{1},X_{2})  \label{eq140} \\
&=&y_{0}\Psi _{00}+y_{1}\Psi _{10}+y_{2}\Psi _{01}+y_{3}\Psi _{11}+y_{4}\Psi
_{20}+y_{5}\Psi _{02},  \nonumber
\end{eqnarray}
where $\Psi _{i_{1}i_{2}}(X_{1},X_{2})=\Psi _{i_{1}}(X_{1})\Psi
_{i_{2}}(X_{2})$ is the product of univariate Hermite polynomials of degree $
i_{1}$ and $i_{2}$ satisfying $i_{1}+i_{2}\leq 2$. Such polynomials can
easily be computed according to the recurrence formula (\ref{eq40}) together
with (\ref{eq50}) with $a_{n}=0$ and $b_{n}=n$. Should we have used Legendre
polynomials instead of Hermite's ones, they would be built by taking $a_{n}=0
$ and $b_{n}=\frac{n^{2}}{(4n^{2}-1)}$. Table \ref{Tab1} shows the bivariate
Hermite polynomial basis constructed from univariate basis $
\{1,X_{1},X_{1}^{2}-1\}$ and $\{1,X_{2},X_{2}^{2}-1\}$ and the link between
the main index $j$ and the multi-index $(i_{1},i_{2})$. For completness, the
same Table also provides bivariate Legendre polynomial basis constructed
from univariate basis $\{1,X_{1},X_{1}^{2}-1/3\}$ and $
\{1,X_{2},X_{2}^{2}-1/3\}$.

\begin{table}[!h]
\begin{tabular}{|c|c|c|c|}
\hline
$j$ & $(i_{1},i_{2})$ & $\Psi _{i_{1}i_{2}}$ (Hermite) & $\Psi _{i_{1}i_{2}}$
(Legendre) \\ \hline
$0$ & $(0,0)$ & $1$ & $1$ \\ \hline
$1$ & $(1,0)$ & $X_{1}$ & $X_{1}$ \\ \hline
$2$ & $(0,1)$ & $X_{2}$ & $X_{2}$ \\ \hline
$3$ & $(1,1)$ & $X_{1}X_{2}$ & $X_{1}X_{2}$ \\ \hline
$4$ & $(2,0)$ & $X_{1}^{2}-1$ & $X_{1}^{2}-1/3$ \\ \hline
$5$ & $(0,2)$ & $X_{2}^{2}-1$ & $X_{2}^{2}-1/3$ \\ \hline
\end{tabular}
\vspace{0.2cm}
\caption{\small Bivariate Hermite and Legendre Polynomial Basis}\label{Tab1}
\end{table}
\end{example}

\subsection{Computing PC coefficients}\label{PC}

$\quad$

The way to compute the PC coefficients $\{y_{j}\}_{0\leq j\leq P_{d}^{p}}$
appearing into equation (\ref{eq130}) can be casted into two different
families: projection methods and regression methods \cite{Pettersson}. Here
we have used the first one which consists in premultiplying (\ref{eq130}) by 
$\Psi _{j}(\mathbf{X})$ and by taking the expectation of the resulting
product. Using the orthogonality of the PC basis, most of the terms cancel
and we end up with 
\begin{equation}
y_{j}=E(Y(\mathbf{X})\Psi _{j}(\mathbf{X}))=\int\limits_{I\subset \mathbb{R}
^{d}}Y(\mathbf{x})\Psi _{j}(\mathbf{x})f(\mathbf{x})d\mathbf{x}\,\,\text{
for } j=0,1..,P_{d}^{p}.  \label{eq150}
\end{equation}

The above integral can be evaluated through different techniques: going from
rough Monte-Carlo sampling simulation to Gaussian quadrature rule or sparse
quadrature rules when the dimension $d$ of the input random is high. Here we
evaluate such integrals using Gaussian quadrature rules which take the form
\begin{equation}
\int\limits_{I\subset \mathbb{R} ^{d}}Y(\mathbf{x})f(\mathbf{x})d\mathbf{x}
\text{ }\approx \sum_{i_{1}=1}^{N_{g_{1}}}\sum_{i_{2}=1}^{N_{g_{2}}}...
\sum_{i_{d}=1}^{N_{g_{d}}}\omega _{i_{1}}\omega _{i_{2}}..\omega _{i_{d}}Y( 
\widetilde{x}_{i_{1}},\widetilde{x}_{i_{2}}..,\widetilde{x}_{i_{d}}),
\label{eq160}
\end{equation}
where $\{\omega _{i_{k}}\}_{1\leq i_{k}\leq N_{g_{k}}}$ are the quadrature
weight and $\{\widetilde{x}_{i_{k}}\}_{1\leq i_{k}\leq N_{g_{k}}}$ are the
quadrature points. The Gaussian quadrature rules are such that the
evaluation of the integral is exact if $Y$ is a multivariate polynomial
containing monomials $x_{i_{k}}$ of maximum degree $2N_{g_{k}}-1$. From
equations (\ref{eq150}) and (\ref{eq160}), we see that it suffices to
evaluate the response $Y$ at $N_{g_{1}}\times N_{g_{2}}...\times N_{g_{\ell}}$
deterministic quadrature points in order to compute the PC coefficients.
Furthermore since $\Psi _{0}=1$, setting $j=0$ in equation (\ref{eq150})
leads to 
\begin{equation}
y_{0}=E(Y(\mathbf{X})),  \label{eq170}
\end{equation}
i.e. the first coefficient of the PC expansion is the expectation of the
random response of the system. Similarly, by considering the approximation
of $Y^{2}$ 
\begin{equation}
Y(\mathbf{X})^{2}\approx
\sum\limits_{i=0}^{P_{d}^{p}}\sum\limits_{j=0}^{P_{d}^{p}}y_{i}y_{j}\Psi
_{i}(\mathbf{X})\Psi _{j}(\mathbf{X}),  \label{eq180}
\end{equation}
and taking the expectation on each side, the orthogonality of the PC basis
leads to 
\begin{equation}
E(Y(\mathbf{X})^{2})=\sum\limits_{i=0}^{P_{d}^{p}}y_{i}^{2},  \label{eq190}
\end{equation}
from which the variance of the random response of the system can easily be
deduced. Furthermore, the PC decomposition provides a convenient way of
computing Sobol' indices, as explained in the next section.

\section{Sensitivity Analysis}\label{sec3}

The purpose of sensitivity analysis, is to investigate the influence of each
input parameter and their possible interactions onto the output measures.
They can be casted into two main families: local analysis based on a local
perturbation around an average value and global analysis that consider input
parameters as random variables and decompose the output variance into
several components. The Sobol' indices belong this last type of family.

\subsection{Sobol' indices}

$\quad$

As in the previous section, we consider the mathematical model $Y=Y(\mathbf{X
})$, where the input parameter $\mathbf{X}=(X_{1},\ldots ,X_{d})$ are $d$
independant random variables belonging to $L_{2}(\Omega _{i},\mathcal{A}
_{i}, \mathcal{P}_{i})$ for $i=1,..,d$ and similarly, $Y$ is assumed to
belong to $L_{2}(\Omega ,\mathcal{A},\mathcal{P})$. In 1993, Sobol \cite
{sobol} proposed an indicator of the influence of the input parameter $X_{i}$
defined by 
\begin{equation}
S_{i}=\frac{V(E(Y/X_{i}))}{V(Y)}=\frac{V_{i}}{V},  \label{eq200}
\end{equation}
commonly termed `first order Sobol' indices'. $V_{i}=V(E(Y/X_{i}))$ is the
conditionnal variance of $Y$ with respect to $X_{i}$ and $V=V(Y)$ is the
total variance of $Y$. Similarly, sensitivity indices of higher order can be
defined by first introducing the following decomposition of the total
variance, which is valid for any output $Y$ belonging to $L_{2}(\Omega ,
\mathcal{A},\mathcal{P})$ 
\begin{equation}
V=\sum_{i_{1}=1}^{d}V_{i}+\sum_{1\leq i_{1}<i_{2}\leq
d}V_{i_{1}i_{2}}+\sum_{1\leq i_{1}<i_{2}<i_{3}\leq
d}V_{i_{1}i_{2}i_{3}}+\cdots+V_{i_{1}\ldots i_{d}},  \label{eq210}
\end{equation}
where 
\begin{eqnarray}
V_{i_{1}} &=&V(E(Y/X_{i_{1}}))  \notag \\
V_{i_{1}i_{2}} &=&V(E(Y/X_{i_{1}},X_{i_{2}}))-V_{i_{1}}-V_{i_{2}}  \notag \\
V_{i_{1}i_{2}i_{3}}
&=&V(E(Y/X_{i_{1}},X_{i_{2}},X_{i_{3}}))-V_{i_{1}i_{2}}-V_{i_{1}i_{3}}-V_{i_{2}i_{3}}-V_{i_{1}}-V_{i_{2}}-V_{i_{3}}
\notag \\
&&\cdots  \label{eq220} \\
V_{i_{1}\ldots i_{d}} &=&V-\sum_{i=1}^{d}V_{i}-\sum_{1\leq i_{1}<i_{2}\leq
d}V_{i_{1}i_{2}}-\cdots-\sum_{1\leq i_{1}<i_{2}\ldots <i_{d-1}\leq
d}V_{i_{1}\ldots i_{d-1}}  \notag
\end{eqnarray}

Then, Sobol' indices of order $k$ are given by 
\begin{equation}
S_{i_{1}\ldots i_{k}}=\frac{V_{i_{1}\ldots i_{k}}}{V}.  \label{eq230}
\end{equation}

Although Sobol' indices could be computed by estimating integrals apearing
into equations (\ref{eq220}), such a procedure would be both computationaly
expensive and hardly tractable. Instead, computing those indices from the PC
representation of the random output $Y$ turns out to be an efficient
alternative and has been the method of choice for many years since the
seminal work of Sudret \cite{sudret}. For that, we define by $
I_{i_{1},i_{2},\ldots,i_{s}}$ ($s\le d$) the set of $d$-dimensional vectors $\boldsymbol{\alpha}=(\alpha
_{1},\ldots,\alpha _{d})$ with $\boldsymbol{\alpha}\not=\boldsymbol{0}$ and $0\le \alpha_{1}\le\ldots\le\alpha _{d}\le d$  that selects elements of the PC basis 
$\Psi _{\alpha _{1}\ldots\alpha _{d}}$
defined by (\ref{eq70}) containing solely the variables $X_{i_1}, X_{i_2}, \ldots,  X_{i_s}$.

This way, the multi-indices defined by $
I_{i}$ will select elements of the PC basis depending only on the variable $
X_{i}$. Similarly, $I_{i,j}$ will select elements of the PC basis that
depend on $X_{i}$ and $X_{j},$ to the exclusion of any other variable. With
this notation, Sobol's indices of order $k$ are simply given as a function
of the PC coefficients as follows: 
\begin{equation}
S_{i_{1}\ldots i_{k}}=\frac{1}{V}\sum_{(\alpha_{1},\ldots,\alpha_{d})\in
I_{i_{1},\ldots,i_{k}}}y_{\alpha_{1},\ldots,\alpha_{d}}^{2}. \label{eq240}
\end{equation}

For a problem with $d$ input random parameters, it can be shown that $
(2^{d}-1)$ Sobol' indices can be computed for each output random quantity
of interest. When the number of input r.v. is large, the number of Sobol'
indices grows exponentially and it becomes difficult to draw information
from these statistics. This is why, in 1996, Homma and Saltelli \cite{Homma}
introduced the total sensitivity indices $S_{T_{i}}$ $(i=1,\ldots,d)$ which
measures the total effect of the $i^{th}$ random input parameter. It is
defined as the sum of all sensitivity indices $S_{i_{1}\ldots i_{k}}$ $
(k=1,\ldots,d)$ for which, one of the indices $i_{1},i_{2},\ldots,i_{k}$ is equal to 
$i$, i.e.
\[
S_{T_i}=\sum_{\substack{ k=1 \\ (i_{1},\ldots,i_{k})\in J_{i}^{k}}}^{d}S_{i_{1},\ldots,i_{k}},
\]
where $J_{i}^{k}$ is the set of $k$-dimensionel vectors $(i_{1},\ldots,i_{k})$ with $1\le i_{1}<\ldots<i_{k}\le d$, such that one 
of its components is equal to $i$. We can also compute these indices from the PC coefficients as follows:
\begin{equation*}
S_{T_i}=1-\frac{1}{V}\sum_{(\alpha_{1},\ldots,\alpha_{d})\in
\overline{I}_{i}}y_{\alpha_{1},\ldots,\alpha_{d}}^{2}, 
\end{equation*}
where $\overline{I}_{i}$ is the complementary set of ${I}_{i}$. 

\begin{example}$\quad$

The computation of the sensitivity indices for the example \ref{ex1} where, $p=2$ and $d=2$ gives 
\begin{eqnarray}
S_{1} &=&\frac{1}{V}(y_{1,0}^{2}+y_{2,0}^{2}),  \nonumber \\
S_{2} &=&\frac{1}{V}(y_{0,1}^{2}+y_{0,2}^{2}),  \label{eq250} \\
S_{1,2} &=&\frac{1}{V}y_{1,1}^{2},  \nonumber
\end{eqnarray}
with 
\begin{equation}
V=y_{1,0}^{2}+y_{0,1}^{2}+y_{1,1}^{2}+y_{2,0}^{2}+y_{0,2}^{2}  \label{eq260}
\end{equation}
In accordance with (\ref{eq210}), we see that the variance $V$ of the output
can be divided into three parts: $V=S_{1}+S_{2}+S_{1,2}.$ The first part $
S_{1}$ represents the influence of the first input random variable alone;
the second $S_{2}$ represents the effect of the second one alone whereas $
S_{1,2}$ accounts for the combined effect of the two input random variables.

Similarly, according to the definition of the total sensitivity indices, we
have
\[
S_{T_{1}}=S_{1}+S_{1,2},
\]
and
\[
S_{T_{2}}=S_{2}+S_{1,2}.
\]
\end{example}


\section{Uncertainty analysis in unreliable queueing models}

\label{sec4} 

In this section, we consider the PC expansion for propagating the
uncertainty in performance measures of queueing models with breakdowns and
repairs, due to epistemic uncertainties in the model input parameters. Section \ref{sec:MG1N} studies 
the M/G/1/N queue with breakdowns,
when the perturbation of a single parameter is introduced. Then, section \ref{subsecm} 
addresses the case of the M/M/1/N queue with breakdowns, for multiple random input parameters.


\subsection{The M/G/1/N queue with breakdowns and repairs}

\label{sec:MG1N}$\quad$ 

Consider a finite capacity M/G/1/N queue with server subject to breakdowns
and repairs. Assume that the customers arrive at the system according to a
Poisson stream with rate $\lambda$, the service times are general and
independent identically distributed with mean $1 / \mu $ and we denote the
service time distribution by $F (x) $. There can be at most $N$ customers in the queue 
(including the one being served), and customers
attempting to enter the queue when there are already $N$ customers present
are lost. The service discipline is assumed to be FCFS. At the beginning of
each service, there is a probability $\theta $ that the server breaks down
(and the customer is sent back to the queue) and enters a repair state. The time of repair is exponentially 
distributed with rate $r$. The only instant of time when a server can breaks down is right at the beginning of a service. This system 
is modeled by a Markov chain $X = \{ X_n : n \in \mathbb{N} \}$, embedded at service
completions and completion of a repair and has state space $S = \{0, 1,
\ldots , N - 1 \}$. Let $\Xi= (a_{i,j})_{i, j \in S}$ be the transition
probability matrix of the Markov chain $X$. The matrix $\Xi$ has the following form:
\begin{equation}  \label{eq:p}
\Xi = \left( 
\begin{array}{ccccccc}
a_{0} & a_{1} & a_{2} & a_{3} & \cdots & a_{N-2} & 
1-\sum\limits_{k=0}^{N-2} a_k \\ 
a_{0} & a_{1} & a_{2} & a_{3} & \cdots & a_{N-2} & 
1-\sum\limits_{k=0}^{N-2} a_k \\ 
0 & a_{0} & a_{1} & a_{2} & \cdots & a_{N-3} & 
1-\sum\limits_{k=0}^{N-3} a_k \\ 
0 & 0 & a_{0} & a_{1} & \cdots & a_{N-4} & 
1-\sum\limits_{k=0}^{N-4} a_k \\ 
&  & \ddots &  & \vdots & \vdots & \vdots \\ 
0 & 0 & 0 & 0 & \cdots & a_{0} & 1-a_{0} \\ 
&  &  &  &  &  & 
\end{array}
\right) ,
\end{equation}
where 
\begin{equation}\label{F}
a_{k} = \theta \int_{0}^{\infty}e^{-\lambda x} \frac{(\lambda x)^{k}}{k!
}dF(x) + ( 1 - \theta ) \frac{r}{r+\lambda}\left(\frac{\lambda}{\lambda+r}
\right)^{k}, \quad k = 0 , \ldots , N-2 .
\end{equation}

Note that the Markov chain $X$ is unichain and therefore the stationary distribution exists. Let $\pi $ denote the
stationary distribution of the queue-length process embedded at service
completions and completion of a repair in the M/G/1/N queue with
breakdowns and repairs (see \cite{abass2016} for details). In the sequel, we
consider $\pi $ as a function of the breakdown probability $\theta $, denoted by $\pi (\theta)$. Assume that the probability $\theta $
is estimated from insufficient statistical data, and hence has uncertainty associated
with it. In the next section, we will discuss a functional approach based on PC expansion for computing the epistemic uncertainty in
stationary distribution $\pi (\theta)$, due to epistemic uncertainties in $\theta$.

\subsubsection{Epistemic uncertainty in queueing system}

\label{unc1}$\quad$

In this model, we will assume that the uncertainty only affects one parameter. To analyze the propagation of
epistemic uncertainty, the mean and the variance of the responses will be
computed using the PC expansion.

Let us consider $\boldsymbol{\pi}(\theta)=({\pi}_{0}(\theta), \cdots,{\pi}_{N-1}(\theta))$ be the random response or the model
output of the queueing system. We
consider the following perturbation model on the uncertain parameter:
\begin{equation}
\theta (\omega )=\widehat{\theta }+\sigma _{\theta }~~\varepsilon (\omega ),
\label{stoch:per}
\end{equation}
where $\theta (\omega )$ is a random variable with uniform density on the
interval $[0,1]$. The parameters $\widehat{\theta }$ and $\sigma _{\theta }$ represent
the mean and the standard deviation of $\theta $ and they can be
estimated by statistical methods. $\varepsilon (\omega )$ can be considered
as a random noise inflicted on $\theta $ and it modelizes epistemic
uncertainty. We assume that $\varepsilon (\omega )$ is a random variable with uniform
distribution on $[-1,1]$. Then, the approximation of any component ${\pi}_{\ell}(\theta )$ of the stationnary distribution
by the PC expansions writes
\begin{equation*}
{\pi}_{\ell}(\theta )\approx \sum_{i=0}^{P_{d}^{p}}y_{i}\Psi
_{i}(\theta ),
\end{equation*}
where $\{\Psi _{i}(\theta )\}_{0\leq i\leq P_{d}^{p}}$ form an orthonormal
Legendre polynomial basis, and $y_{i}$ are the coefficients of the
approximating series expansion, see (\ref{eq130}) .

\subsubsection{Numerical results}

$\quad$

In this section, the parameters $N$, $\lambda $ and $r$ are deterministic
and their values are set to $N=7$, $\lambda =1$, $r=0.4$, respectively. The
only input random parameter is $\theta $ and it is given by
\begin{equation*}
\theta =0.5+0.28~~\varepsilon ,~~~~\varepsilon \sim \mathcal{U}(-1,1).
\end{equation*}

For numerical computations, we propose two types of distribution for the service time,
namely, Erlang $(E_{2})$ and Hyperexponential of second order $(H_{2})$.
Those two distributions are described in the following subsections.

\subsubsection{Epistemic uncertainty in the {\rm M/${\rm E}_{2}$/1} queue with server
breakdowns and repairs}

$\quad$

In (\ref{F}), the density $F'(x)=f(x)$ of the service time is assumed to be Erlang of second order, i.e.
\begin{equation*}
f(x)=\frac{\mu _{1}\mu _{2}}{\mu _{2}-\mu _{1}}(e^{-\mu _{1}x}-e^{-\mu _{2}x}){\bf 1}_{[0,+\infty)}(x), 
\end{equation*}
where the rates $\mu _{1}=4$ and $\mu _{2}=2$.

We determine the stationary distribution vector with respect to the random
variable $\theta (\omega )$ by solving the system $\pi _{\theta
}\Xi_{\theta }=\pi _{\theta }$ and $\sum \pi _{\theta }=1$, where $\Xi_{\theta }$
denotes the transition probability given by (\ref{eq:p}). Therefore,
the random outputs of interest are the stationary distributions ${
\pi }_{i}$, $i\in \{0,\ldots,6\}.$ Their projection is performed on a
monovariate PC basis of degree $n=4$ and the coefficients are
computed using a $6$-points Gaussian quadrature rule as explained in section \ref{PC}.
 Numerical results for the expectations and the variances using a
Monte-Carlo simulations of a sample size $N_{MC}=1000$ and $N_{MC}=100000$
are given in Table \ref{exp1} and Table \ref{vari1}, respectively. Results are
compared with those of the PC expansion. We note that the Monte Carlo
simulations (MC in short) converge to the PC results; however the later are obtained at a
fraction of the cost of the MC simulations. Similarly moments of
higher order (skewness and kurtosis) are given in Tables \ref{skew11} and \ref{kur11}, and the same
conclusions can be drawn.

\begin{table}[!h]
\centering
\begin{tabular}{|c|c|c|c|}
\hline
$E({\pi}_{i})$ & PC& MC  {\scriptsize($ N_{MC}
=1000$)} & MC  {\scriptsize($N_{MC}=100 000$)}\\ \hline
${\pi}_{0}$ & $0.0237$ & $0.0230$ & $0.0237$\\ \hline
${\pi}_{1}$ & $0.0413$ & $0.0422$ & $0.0413$ \\ \hline
${\pi}_{2}$ & $0.0660$ & $0.0667$ & $0.0660$ \\ \hline
${\pi}_{3}$ & $0.1018$ & $0.1026$ & $0.1018$ \\ \hline
${\pi}_{4}$ & $0.1557$ & $0.1563$ & $0.1558$ \\ \hline
${\pi}_{5}$ & $0.2393$ & $0.2401$ & $0.2392$ \\ \hline
${\pi}_{6}$ & $0.3722$ & $0.3719$ & $0.3722$ \\ \hline
\end{tabular}
\vspace{0.2cm}
\caption{\small Expected value of the steady-state vector in M/${\rm E}_{2}$/1/7 queue}
\label{exp1}
\end{table}

\begin{table}[!h]
\centering
\begin{tabular}{|c|c|c|c|}
\hline
$V({\pi}_{i})$ & PC$\times 10^{-4}$ & MC$\times 10^{-4}$  {\scriptsize($N_{MC}=1000$)}
& MC$\times 10^{-4}$   {\scriptsize($N_{MC}=100000$)} \\ \hline
${\pi}_{0}$ & $1.7753$ & $1.6686$ & $1.7809$
\\ \hline
${\pi}_{1}$ & $3.6413$ & $3.6789$ & $3.6423$
\\ \hline
${\pi}_{2}$ & $5.2349$ & $5.1333$ & $5.2175$
\\ \hline
${\pi}_{3}$ & $5.0800$ & $4.9922$ & $5.0508$
\\ \hline
${\pi}_{4}$ & $1.8737$ & $1.8585$ & $1.8732$
\\ \hline
${\pi}_{5}$ & $2.1440$ & $2.0101$ & $2.1534$
\\ \hline
${\pi}_{6}$ & $58.7699$  & $61.0693$   & $58.7199$ \\ \hline
\end{tabular}
\vspace{0.2cm}
\caption{\small Variance of the steady-state vector in M/${\rm E}_{2}$/1/7 queue}
\label{vari1}
\end{table}
\begin{table}[!h]
\centering
\begin{tabular}{|c|c|c|c|}
\hline
${\rm Skew}({\pi}_{i})$ & PC& MC  {\scriptsize($N_{MC}=1000$)} & MC  {\scriptsize($N_{MC}=100000$)} \\ \hline
${\pi}_{0}$ & $0.5053$ & $0.5429$ & $0.5064$ \\ \hline
${\pi}_{1}$ & $0.3565$ & $0.2911$ & $0.3550$ \\ \hline
${\pi}_{2}$ & $0.1730$ & $0.1435$ & $0.1701$ \\ \hline
${\pi}_{3}$ & -$0.0722$ & -$0.0927$ & -$0.0730$ \\ \hline
${\pi}_{4}$ & -$0.5848$ & -$0.6616$ & -$0.5922$ \\ \hline
${\pi}_{5}$ & -$0.7140$ & -$0.7968$ & -$0.7036$ \\ \hline
${\pi}_{6}$ & $0.0543$ & $0.0692$ & $0.0540$ \\ \hline
\end{tabular}
\vspace{0.2cm}
\caption{\small Skewness coefficient of the steady-state vector in M/${\rm E}_{2}$/1/7
queue }
\label{skew11}
\end{table}

\begin{table}[!h]
\centering
\begin{tabular}{|c|c|c|c|}
\hline
${\rm Kurt}({\pi}_{i})$ & PC & MC  {\scriptsize($N_{MC}=1000$)} & MC  {\scriptsize($N_{MC}=100000$)} \\ \hline
${\pi}_{0}$ & $2.0384$ & $2.1425$ & $2.0341$ \\ \hline
${\pi}_{1}$ & $1.8968$ & $1.9122$ & $1.8953$ \\ \hline
${\pi}_{2}$ & $1.7864$ & $1.7682$ & $1.7900$ \\ \hline
${\pi}_{3}$ & $1.7463$ & $1.7513$ & $1.7495$ \\ \hline
${\pi}_{4}$ & $2.0301$ & $2.0290$ & $2.0385$ \\ \hline
${\pi}_{5}$ & $2.1921$ & $2.3603$ & $2.1972$ \\ \hline
${\pi}_{6}$ & $1.7759$ & $1.8007$ & $1.7724$ \\ \hline
\end{tabular}
\vspace{0.2cm}
\caption{\small Kurtosis coefficient of the steady-state vector in M/${\rm E}_{2}$/1/7
queue }
\label{kur11}
\end{table}

\begin{figure}[!h]
\begin{center}
\includegraphics[width=15cm,height=10cm]{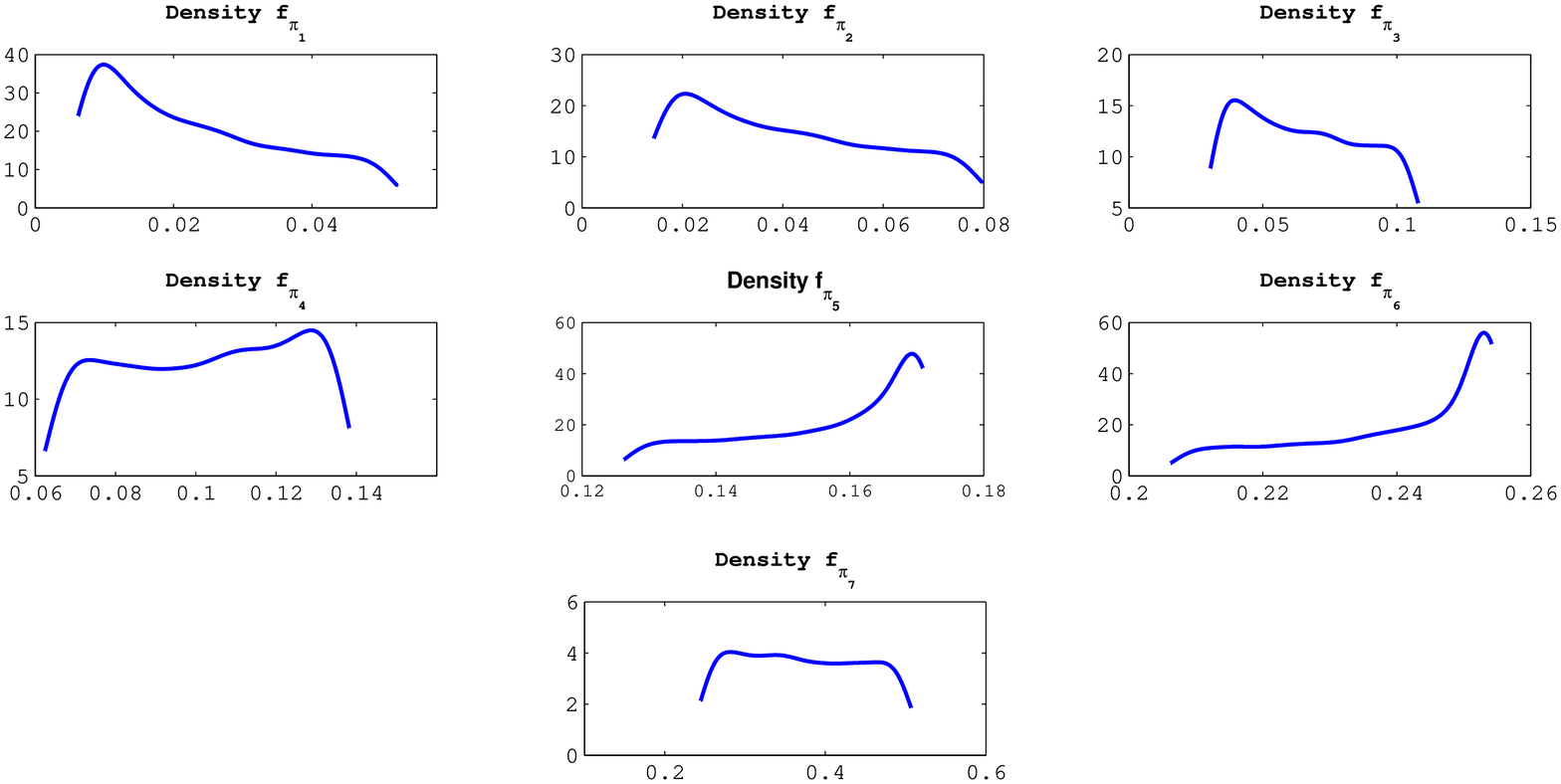}
\end{center}
\caption{\small Mariginal probability density $f_{\pi_i}$ of the steady-state vector $\protect\pi(\theta)$ in M/${\rm E}_{2}$/1/7 queue}
\label{density1}
\end{figure}

Figure \ref{density1} show the density of the output stationary distributions for a
random uniform input of the form (\ref{stoch:per}). The shape of the densities
show that they are far from being uniform. Furthermore, the plots are
coherent with Tables \ref{skew11} and \ref{kur11}. Indeed, a positive skewness reflects a
random variable with higher probabilities in the left part of its support,
and vice versa. 


\subsubsection{Epistemic uncertainty in the {\rm M/${\rm H}_{2}$/1} queue with server breakdowns and repairs}

$\quad$

In this section, the density $F'(x)=f(x)$ of the service time is now assumed to be
Hyperexponential of second order i.e. 
\begin{equation*}
f(x)=(\gamma \mu _{1}e^{-\mu _{1}x}+(1-\gamma )\mu _{2}e^{-\mu _{2}x}){\bf 1}_{[0,\infty)}(x),
\label{eq270}
\end{equation*}
where the rates are $\mu _{1}=3/2$, $\mu _{2}=3$; and $\gamma =0.3$. All
the other parameters are set as in the previous section, including the
random input. Similarly, Tables \ref{exp12}, \ref{vari12}, \ref{skew12} and \ref{kur12}  show the moments of the stationary
distributions when the Hyperexponential law is used instead of the Erlang
one for the service time. Here also, we can observe the good convergence
properties of the PC, showing the robustness of the method independently of
the law of the service time.

\begin{table}[!h]
\centering
\begin{tabular}{|c|c|c|c|}
\hline
$E({\pi}_{i})$ & PC & MC {\scriptsize($N_{MC}=1000$)} & MC  {\scriptsize($N_{MC}=100000$)} \\ \hline
${\pi}_{0}$ & $0.0397$ & $0.0407$ & $0.0396$ \\ \hline
${\pi}_{1}$ & $0.0602$ & $0.0588$ & $0.0602$ \\ \hline
${\pi}_{2}$ & $0.0853$ & $0.0847$ & $0.0853$ \\ \hline
${\pi}_{3}$ & $0.1179$ & $0.1179$ & $0.1179$ \\ \hline
${\pi}_{4}$ & $0.1619$ & $0.1615$ & $0.1619$ \\ \hline
${\pi}_{5}$ & $0.2233$ & $0.2238$ & $0.2233$ \\ \hline
${\pi}_{6}$ & $0.3118$ & $0.3103$ & $0.3117$ \\ \hline
\end{tabular}
\vspace{0.2cm}
\caption{\small Expected value of the steady-state vector in M/${\rm H}_{2}$/1/7 queue}
\label{exp12}
\end{table}

\begin{table}[!h]
\centering
\begin{tabular}{|c|c|c|c|}
\hline
$ V({\pi}_{i})$ & PC $\times 10^{-4}$ & MC$\times 10^{-4}$  {\scriptsize($N_{MC}=1000$)}
& MC$\times 10^{-4}$  {\scriptsize($N_{MC}=100000$)} \\ \hline
${\pi}_{0}$ & $4.7367$ & $4.5451$ & $4.7370$
\\ \hline
${\pi}_{1}$ & $6.7973$ & $6.5517$ & $6.8120$
\\ \hline
${\pi}_{2}$ & $7.1124$ & $7.1342$ & $7.1237$ \\ \hline
${\pi}_{3}$ & $4.7642$ & $4.5830$ & $4.7726$
\\ \hline
${\pi}_{4}$ & $0.7665$ & $0.8177$ & $0.7653$
\\ \hline
${\pi}_{5}$ & $5.4708$ & $5.2979$ & $5.4699$
\\ \hline
${\pi}_{6}$ & $64.0300$ & $62.3300$ & $64.0396$ \\ \hline
\end{tabular}
\vspace{0.2cm}
\caption{\small Variance of the steady-state vector in M/${\rm H}_{2}$/1/7 queue}
\label{vari12}
\end{table}

\begin{table}[!h]
\centering
\begin{tabular}{|c|c|c|c|}
\hline
${\rm Skew}({\pi}_{i})$ & PC & MC  {\scriptsize($N_{MC}=1000$)} & MC  {\scriptsize($N_{MC}=100000$)} \\ \hline
${\pi}_{0}$ & $0.3481$ & $0.2807$ & $0.3531$ \\ \hline
${\pi}_{1}$ & $0.1478$ & $0.2379$ & $0.1474$ \\ \hline
${\pi}_{2}$ & -$0.0881$ & -$0.0925$ & -$0.0896$ \\ \hline
${\pi}_{3}$ & -$0.4210$ & -$0.4300$ & -$0.4230$ \\ \hline
${\pi}_{4}$ & -$1.2389$ & -$1.1806$ & -$1.2381$ \\ \hline
${\pi}_{5}$ & -$0.3119$ & -$0.3561$ & -$0.3101$ \\ \hline
${\pi}_{6}$ & $0.2751$  & $0.3242$ & $0.2747$ \\ \hline
\end{tabular}
\vspace{0.2cm}
\caption{\small Skewness coefficient of the steady-state vector in M/${\rm H}_{2}$/1/7
queue}
\label{skew12}
\end{table}
\begin{table}[!h]
\centering
\begin{tabular}{|c|c|c|c|}
\hline
${\rm Kurt}({\pi}_{i})$ & PC& MC   {\scriptsize($N_{MC}=1000$)} & MC  {\scriptsize($N_{MC}=100000$)} \\ \hline
${\pi}_{0}$ & $1.8643$ & $1.9138$ & $1.8624$ \\ \hline
${\pi}_{1}$ & $1.7524$ & $1.7200$ & $1.7568$ \\ \hline
${\pi}_{2}$ & $1.7297$ & $1.6881$ & $1.7266$ \\ \hline
${\pi}_{3}$ & $1.8885$ & $1.8168$ & $1.8935$ \\ \hline
${\pi}_{4}$ & $3.3357$ & $3.7223$ & $3.3412$ \\ \hline
${\pi}_{5}$ & $1.7598$ & $1.7821$ & $1.7586$ \\ \hline
${\pi}_{6}$ & $1.8392$ & $1.8446$ & $1.8372$ \\ \hline
\end{tabular}
\vspace{0.2cm}
\caption{\small Kurtosis coefficient of the steady-state vector in M/${\rm H}_{2}$/1/7
queue}
\label{kur12}
\end{table}

\begin{figure}[!h]
\begin{center}
\includegraphics[width=15cm,height=10cm]{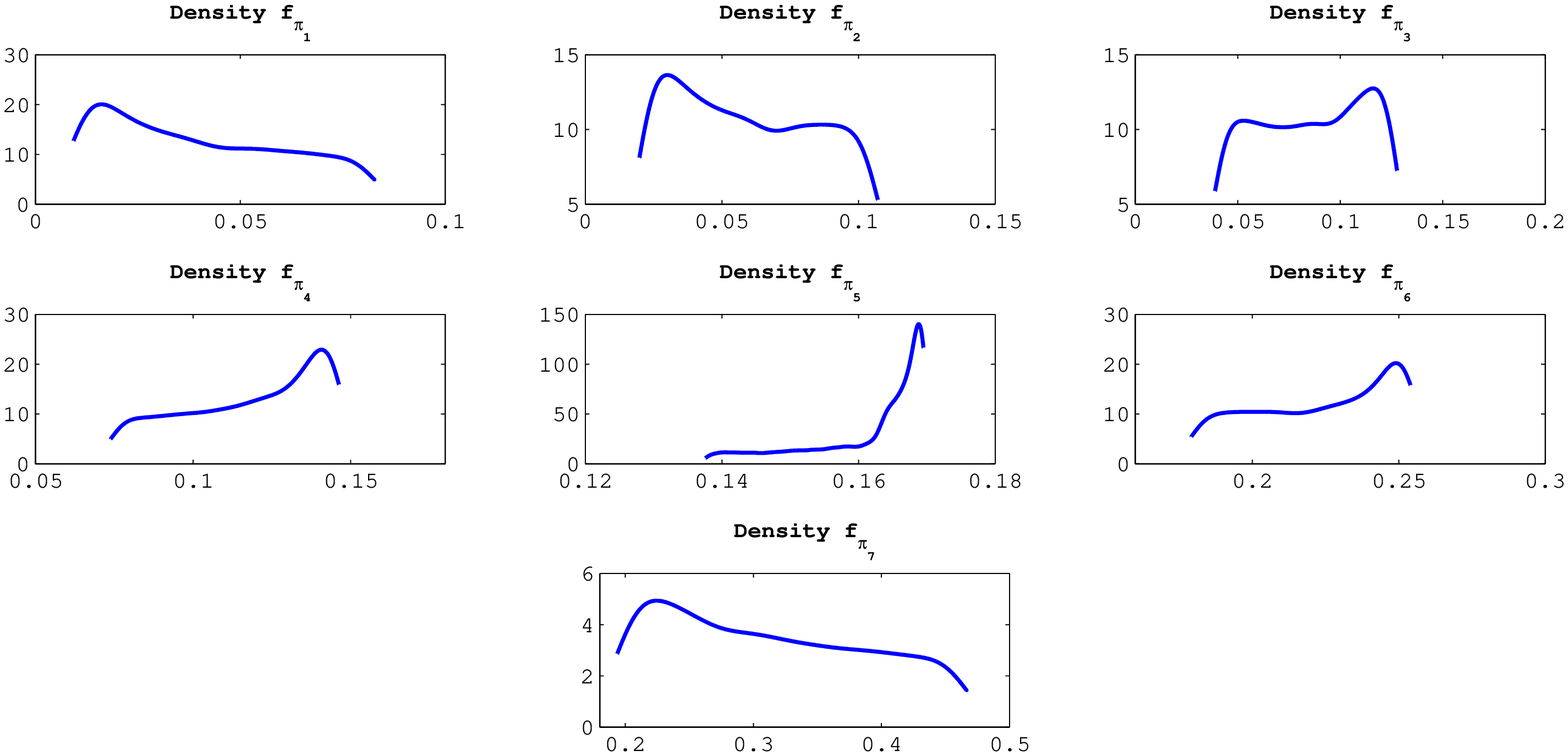}
\end{center}
\caption{\small Mariginal probability density $f_{\pi_i}$ of of the steady-state vector in M/${\rm H}_{2}$/1/7 queue}
\end{figure}

\subsection{The M/M/1 queue with server breakdowns and threshold-based
recovery policy}

\label{subsecm}$\quad$

Consider a finite capacity M/M/1 queue with server breakdowns and
threshold-based recovery policy. If a server provides service for a
customer, the server may experience a breakdown, and the repair begins when
the number of customers present in the system exceeds some prespecified
threshold level $q$ ($1\leq q\leq N$). The stream customers arrive at the
queue according to a Poisson law of parameter $\lambda $. Arriving customers
form a single waiting line based on the order of their arrivals and the server
can serve only one customer at a time. There can be at most $N$ customers
present at the queue (including the one in service), and customers
attempting to enter the queue when $N$ customers are already present are
lost. A single server with exponential distributed service times with rate $
\mu $ is considered. We assume that the server can break down only if the
system is not empty. The lengths of breakdowns are identically distributed
and follow an exponential distribution with rate $\alpha $. The duration of
reparation of the server is assumed to follow an exponential law with
parameter $\beta $. After the server being repaired, it switches to working
state and continues to provide service until the system becomes empty. The
arrival flow customers, service times, breakdown times and repair times are
assumed to be mutually independent input random variables. The service
discipline is FCFS.

A typical state of this system may be denoted by $X(t)=\{(N(t),Y(t));t\geq 0\}$
, where $N(t)$ is the number of customers in the queue at time $t$, and $Y(t)
$ is a random variable representing the server state at time $t$. If at time 
$t$, the server is experiencing a breakdown, then $Y(t)=1$; otherwise, the
server is in the working state and $Y(t)=0$. The stochastic process $X(t)$
is a continuous time Markov chain whose state space $S=\{(n,0):n=0,\ldots
,N\}\cup \{(n,1):n=1,\ldots ,N\}$. The infinitesimal generator matrix $Q$ of
the continuous time Markov chain $X(t)$ has the following block-tridiagonal
form: 
\begin{equation}
Q=\left( 
\begin{array}{cccccccccc}
A_{0} & B_{0} & 0 & 0 & 0 & 0 & 0 & 0 & 0 & \ldots    \\ 
D_{1} & {B}_{1} & {F}_{1} & 0 & 0 & 0 & 0 & 0 & 0 & \ldots    \\ 
0 & {D}_{2} & {B}_{2} & {F}_{2} & 0 & 0 & 0 & 0 & 0 & 
\ldots   \\ 
0 & 0 & {D}_{3} & {B}_{3} & {F}_{3} & 0 & 0 & 0 & 0 & 
\ldots    \\ 
\ddots  & \ddots  & \ddots  & \ddots  & \ddots  & \ddots  & \ddots  & \ddots 
& \ddots  & \ddots   \\ 
0 & 0 & 0 & \ldots  & {D}_{s} & M_{s} & {F}_{s} & 0 & 0 & \ldots   \\ 
0 & 0 & 0 & 0 & \ldots  & {D}_{s+1} & M_{s+1} & {F}_{s+1} & 0 & 
\ldots    \\ 
\ddots  & \ddots  & \ddots  & \ddots  & \ddots  & \ddots  & \ddots  & \ddots 
& \ddots  & \ddots    \\ 
0 & 0 & 0 & 0 & 0 & 0 & 0 & \ldots  & {D}_{N} & E_{N}   \\ 
\end{array}
\right),
\end{equation}
where $B_{0}=\lambda,  A_{0} =-\lambda, D_{1} = \left( \begin{array}{c}\mu  \\ 0  \\ \end{array}\right), {B}_{i} = \left( 
\begin{array}{cc}
-(\lambda+\mu+\alpha) & \alpha \\ 
0 & -\lambda \\  
\end{array}
\right),~~~~ \,i=1,\cdots q-1,$

\begin{equation*}
 {F}_{j}= \left( 
\begin{array}{cc}
\lambda & 0 \\ 
0 & \lambda \\  
\end{array}
\right),\, \, j=1,\cdots,N-1,\quad {D}_{k}= \left( 
\begin{array}{cc}
\mu & 0 \\ 
0 & 0 
\end{array}
\right), \, k=2,\cdots N
\end{equation*}

\begin{equation*}
 M_{s}= \left( 
\begin{array}{cc}
-(\lambda+\alpha+\mu) & \alpha \\ 
\beta & -(\beta+\lambda) \\ 
\end{array}
\right), ~~~~\,s=q,\cdots N-1,\quad
E_{N}= \left( 
\begin{array}{cc}
-(\alpha+\mu) & \alpha \\ 
\beta & -\beta \\ 
\end{array}
\right).
\end{equation*}

The Markov transition diagram of the process is given below
\begin{figure}[!h]
\begin{center}
\includegraphics[width=14cm,height=5cm]{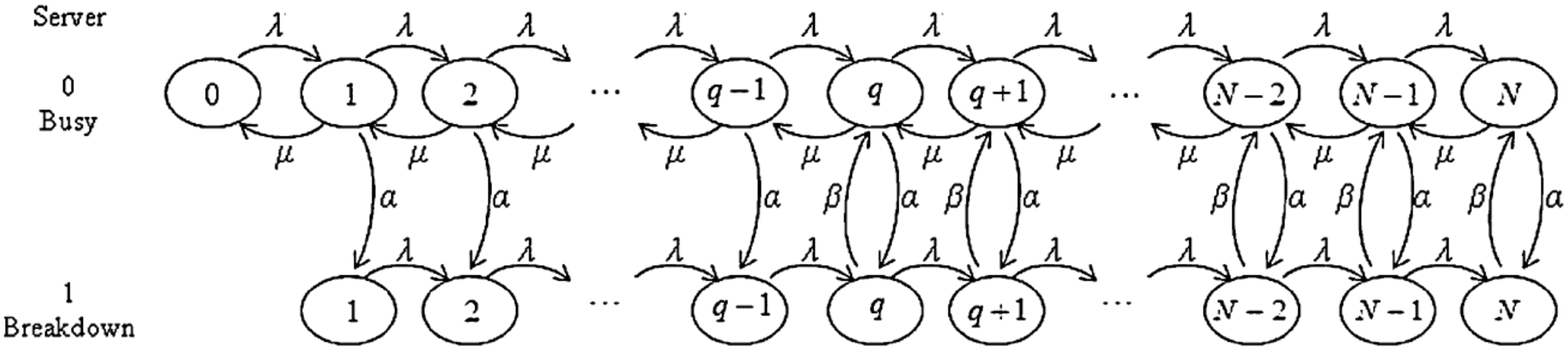}
\end{center}
\vspace{-1cm}
\caption{Markov transition diagram of the M/M/1/N queue with server breakdowns and threshold-based recovery policy}
\end{figure}

Denote by $\pi _{m,n}$, $m=0,1,~~n\geq 0$ the stationary distribution of
the Markov chain $X(t)$. We have $\pi _{0}=\pi _{0,0}$, and $\pi _{n}=(\pi
_{0,n},\pi _{1,n}),n=1,\ldots ,N$. The steady-state analysis of the
Markov chain $X(t)$ can be performed via a recursive scheme \cite{Yang:2013}
, or by using matrix analytic methods \cite{Neuts}. In the sequel, we follow
a different method, and we will provide a functional approach based on the PC expansion 
for obtaining the stationary distribution $\pi
_{m,n}$ in terms of some parameters, which are obtained under the
epistemic uncertainties.

\subsubsection{Epistemic uncertainty in queueing system}\label{unc}$\quad$

In this model, we consider the epistemic uncertainty in all the input
parameters of the queueing system, $\lambda $, $\mu $, $\alpha $, $\beta $.
However, now $\pi _{m,n}(\boldsymbol{\theta })$ is assumed to be the random
output response of the model of the queueing system, where  $\boldsymbol{\theta }
=[\theta _{1}=\lambda ,\theta _{2}=\mu ,\theta _{3}=\alpha ,\theta
_{4}=\beta ]$ is the vector that describes the model parameters. We consider
an epistemic uncertainty in the input parameter $\boldsymbol{\theta }$ of known
probability density function. Due to the uncertainty in input parameters,
the responses $\pi _{m,n}(\boldsymbol{\theta }(\omega ))$ are considered as
random variables. The uncertainty analysis of the functional $\pi _{m,n}$
may be quantified by the computation of its mean and variance and moments of
higher order. 
\medskip

Furthermore, we proceed to a sensitivity analysis in order to
determine the most sensitive parameters with respect to the performance
measures. A convenient way to compute those quantities is to approximate the
functional $\pi _{m,n}(\boldsymbol{\theta }(\omega ))$ in the form of the PC
expansion. The sensitivity analysis introduced in section \ref{sec3} help us
to point out which parameters generate a higher sensitivity with respect to the stationary distribution. We
consider the following perturbation models for all the uncertain parameters:
\begin{equation*}
\theta _{j}(\omega )=\hat{\theta _{j}}+\sigma _{j}\varepsilon _{j}(\omega
),\qquad j=1,\cdots,4  \label{stoch:per2}
\end{equation*}
where $\theta _{j}(\omega )$ is random variable, $\hat{\theta _{j}}$, $
\sigma _{j}$, are respectively, estimated mean and the standard deviation of
the random variable $\theta _{j}$ obtained by statistical methods. $
\varepsilon _{j}(\omega )$ is a random noise inflicted on $\theta _{i}$
which models epistemic uncertainty. Here, $\varepsilon _{j}(\omega )$ is assumed
to follow a standard normal distribution. Let $\boldsymbol{\theta }=[\theta
_{1}(\omega )=\lambda ,\theta _{2}(\omega )=\mu ,\theta _{3}(\omega )=\alpha
,\theta _{4}(\omega )=\beta ]$, $\pi _{m,n}(\boldsymbol{\theta }(\omega ))$ and $
\varepsilon =[\epsilon _{1},\ldots ,\epsilon _{4}]$ be random fields.

The approximation of the measure $\pi _{m,n}(\boldsymbol{\theta }(\omega ))$ in
the form of multivariate PC expansions is given by:
\begin{equation*}
\pi _{m,n}(\boldsymbol{\theta })\approx \sum_{i=0}^{P_{4}^{p}}y_{i}\Psi
_{i}(\boldsymbol{\theta }),
\end{equation*}
where $\Psi _{i}(\boldsymbol{\theta })$ is the multivariate orthonormal Hermite
polynomial and $y_{i}$ are the coefficients of the approximating series
expansion, (\ref{eq130}).

\subsubsection{Numerical application}$\quad$

As an illustrative example, the queueing system described in section \ref{subsecm} is
studied. The capacity is set to $N=5$ and the threshold is equal to $q=3$.
In order to study the sensitivity of the queue, we assume that the exact
value of all the input parameters, $\alpha $, $\beta $,$\lambda $, $\mu $,
are not precisely known and we propose the following perturbation on these
parameters:
\begin{equation*}
\alpha=\hat{\alpha} + \sigma_{\alpha} \varepsilon_{1}, ~~
\varepsilon_{1}\sim \mathcal{N}(0,1)\qquad\qquad
\beta=\hat{\beta} + \sigma_{\beta} \varepsilon_{2}, ~~ \varepsilon_{2}\sim 
\mathcal{N}(0,1)
\end{equation*}
\begin{equation*}
\lambda=\hat{\lambda} + \sigma_{\lambda} \varepsilon_{3}, ~~
\varepsilon_{3}\sim \mathcal{N}(0,1)\qquad\qquad
\mu=\hat{\mu} + \sigma_{\mu} \varepsilon_{4}, ~~ \varepsilon_{4}\sim 
\mathcal{N}(0,1)
\end{equation*}

so $\alpha $, $\beta $, $\lambda $, $\mu $, are Gaussian random
variables. The mean and the standard deviation of these random variables are
set respectively to $\hat{\mu}=7.3$, $\hat{\lambda}=2$, $\hat{\alpha}=3$, $
\hat{\beta}=4$, $\sigma _{\alpha }=0.04$, $\sigma _{\beta }=0.02$, $\sigma
_{\lambda }=0.04$ and $\sigma _{\mu }=0.02$.


\begin{table}[!h]
\centering
\begin{tabular}{|c|c|c|c|c|}
\hline
$\pi_{m,n}\backslash S_{\theta_{i}}$ & $S_{\alpha}$ & $S_{\mu}$ & $S_{\beta}$ & $
S_{\lambda}$ \\ \hline
$\pi_{0,0}$ & $0.3799$ & $0.0366$ & $0.0046$ & $0.5789$ \\ \hline
$\pi_{0,1}$ & $0.7942$ & $0.0047$ & $0.0096$ & $0.1913$ \\ \hline
$\pi_{1,1}$ & $0.0144$ & $0.0038$ & $0.0078$ & $0.9740$ \\ \hline
$\pi_{0,2}$ & $0.0632$ & $0.0160$ & $0.0080$ & $0.9126$ \\ \hline
$\pi_{1,2}$ & $0.1650$ & $0.0003$ & $0.0109$ & $0.8236$ \\ \hline
$\pi_{0,3}$ & $0.0445$ & $0.0473$ & $0.0073$ & $0.9008$ \\ \hline
$\pi_{1,3}$ & $0.7812$ & $0.0151$ & $0.0418$ & $0.1619$ \\ \hline
$\pi_{0,4}$ & $0.0502$ & $0.0264$ & $0.0006$ & $0.9229$ \\ \hline
$\pi_{1,4}$ & $0.3120$ & $0.0124$ & $0.0411$ & $0.6345$ \\ \hline
$\pi_{0,5}$ & $0.0528$ & $0.0198$ & $0.0044$ & $0.9229$ \\ \hline
$\pi_{1,5}$ & $0.1512$ & $0.0085$ & $0.0396$ & $0.8005$ \\ \hline
\end{tabular}
\vspace{0.2cm}
\caption{\small First order Sobol' indices for stationary distribution vector in M/M/1/5 queue 
with server breakdowns and threshold-based recovery policy }
\label{so1}
\end{table}


\begin{table}[!h]
\centering
\begin{tabular}{|c|c|c|c|c|c|c|}
\hline
$\pi_{m,n} \backslash S_{\theta_{i},\theta_{j}}$ & $S_{\alpha,\mu}\times10^{-8}$ & $
S_{\alpha,\beta}\times10^{-8}$ & $S_{\alpha,\lambda}\times10^{-6}$ & $S_{\mu,\beta}\times10^{-9}$ & $S_{\mu,\lambda}\times10^{-8}$ & $S_{\beta,\lambda}\times10^{-7}$ \\ 
\hline
$\pi_{0,0}$ & $0.6972$ & $1.3636$ & $0.2051$ & $1.8411$ & $4.9130$ & $1.1427$
\\ \hline
$\pi_{0,1}$ & $0.0144$ & $2.8434$ & $0.0127$ & $0.4145$ & $0.0016$ & $0.6002$
\\ \hline
$\pi_{1,1}$ & $0.0361$ & $0.0175$ & $0.5309$ & $0.3338$ & $0.0569$ & $1.9217$
\\ \hline
$\pi_{0,2}$ & $0.0516$ & $0.7513$ & $0.0183$ & $0.0181$ & $0.0249$ & $0.8644$
\\ \hline
$\pi_{1,2}$ & $0.0263$ & $0.0365$ & $0.0143$ & $0.0109$ & $0.0013$ & $7.3456$
\\ \hline
$\pi_{0,3}$ & $0.8840$ & $0.0198$ & $0.3924$ & $0.0283$ & $0.4031$ & $0.8201$
\\ \hline
$\pi_{1,3}$ & $9.8593$ & $1.6277$ & $0.3311$ & $0.1568$ & $0.0997$ & $0.0147$
\\ \hline
$\pi_{0,4}$ & $0.1210$ & $7.1050$ & $4.0803$ & $4.9641$ & $0.0012$ & $8.0038$
\\ \hline
$\pi_{1,4}$ & $0.0137$ & $0.0368$ & $0.4395$ & $0.9306$ & $0.6662$ & $8.6512$
\\ \hline
$\pi_{0,5}$ & $0.0113$ & $0.2362$ & $0.5077$ & $0.6212$ & $0.0029$ & $0.1801$
\\ \hline
$\pi_{1,5}$ & $0.0250$ & $0.0959$ & $0.0187$ & $0.0433$ & $0.0010$ & $0.0427$
\\ \hline
\end{tabular}
\vspace{0.2cm}
\caption{\small Second order Sobol' indices for stationary distribution vector in M/M/1/5 queue 
with server breakdowns and threshold-based recovery policy}
\label{so3}
\end{table}

\begin{table}[!h]
\centering
\begin{tabular}{|c|c|c|c|c|}
\hline
$\pi_{m,n} \backslash S_{\theta_{i}}^{T}$ & $S_{\alpha}^{T}$ & $S_{\mu}^{T}$ & $S_{\beta}^{T}$ & $S_{\lambda}^{T}$ \\ \hline
$\pi_{0,0}$ & $0.3799$ & $0.0366$ & $0.0046$ & $0.5789$ \\ \hline
$\pi_{0,1}$ & $0.7943$ & $0.0047$ & $0.0096$ & $0.1914$ \\ \hline
$\pi_{1,1}$ & $0.0144$ & $0.0038$ & $0.0078$ & $0.9741$ \\ \hline
$\pi_{0,2}$ & $0.0634$ & $0.0160$ & $0.0080$ & $0.9128$ \\ \hline
$\pi_{1,2}$ & $0.1652$ & $0.0002$ & $0.0109$ & $0.8237$ \\ \hline
$\pi_{0,3}$ & $0.0445$ & $0.0473$ & $0.0074$ & $0.9009$ \\ \hline
$\pi_{1,3}$ & $0.7812$ & $0.0151$ & $0.0419$ & $0.1619$ \\ \hline
$\pi_{0,4}$ & $0.0502$ & $0.0264$ & $0.0005$ & $0.9229$ \\ \hline
$\pi_{1,4}$ & $0.3121$ & $0.0124$ & $0.0411$ & $0.6345$ \\ \hline
$\pi_{0,5}$ & $0.0528$ & $0.0198$ & $0.0044$ & $0.9230$ \\ \hline
$\pi_{1,5}$ & $0.1514$ & $0.0085$ & $0.0396$ & $0.8007$ \\ \hline
\end{tabular}
\vspace{0.2cm}
\caption{\small Total order Sobol' indices for steady-state vector in M/M/1/5 queue 
with server breakdowns and threshold-based recovery policy}
\label{so2}
\end{table}

\bigskip

Table \ref{so1} gives the first order Sobol' indices for all the components of
the stationary distribution vector. From this table, we see that for all but 
$\pi _{0,1}$ and $\pi _{1,3}$, the most influencial input random
parameter is $\lambda $. Table \ref{so3}, which gives the second order Sobol'
indices, shows that the combined efffect of two input random parameter is
negligible in comparison with the effect of a single parameter. The same is
also true for higher order indices since the total Sobol' indices are
almost equal to first order ones, as shown by Table \ref{so2}. Note that for the
components $\pi _{0,1}$ and $\pi _{1,3},$ it is the parameter $\alpha $
which is most influencial.
\medskip

When considering $\pi _{0,1}$ and $\pi _{1,3}$, we can set $\alpha $ as
being the only random parameters and when considering all the other
components of the stationary distribution, we can consider the parameter $
\lambda $ as sole input random variable. This is what we do as a second
numerical experiments. Then, Tables \ref{exp} to \ref{kur2} compare the statistical
moments of the stationary distribution vector when the four input parameters
are random (first column) to the case where only the most influencial
parameter is random (second column). As expected, the two columns of the
tables match rather well.


\begin{table}[!h]
\centering
\begin{tabular}{|c|c|c|c|}
\hline
$E(\pi_{m,n})$  & PC (1 r.v.)&   PC (4 r.v.) \\
\hline
$\pi_{0,0}$ &  $0.2639$ & $0.2639$   \\\hline
$\pi_{0,1}$ &  $0.0723$ & $0.0723$  \\\hline
$\pi_{1,1}$ & $0.1084$ & $0.1084$ \\\hline
$\pi_{0,2}$ & $0.0495$ & $0.0495$ \\\hline
$\pi_{1,2}$ &  $0.1827$  & $0.1827$  \\\hline
$\pi_{0,3}$ &  $0.0636$  & $0.0636$ \\\hline
$\pi_{1,3}$ &  $0.0927$ & $0.0927$  \\\hline
$\pi_{0,4}$ & $0.0428$ &  $0.0428$ \\\hline
$\pi_{1,4}$ &  $0.0523$ & $0.0523$  \\\hline
$\pi_{0,5}$ &  $0.0261$ & $0.0261$ \\\hline
$\pi_{1,5}$ &  $0.0457$ & $0.0457$ \\\hline
\end{tabular}
\vspace{0.2cm}
\caption{\small Expected Value of the stationary distribution in M/M/1/5 queue 
with server breakdowns and threshold-based recovery policy}
\label{exp}
\end{table}

\begin{table}[!h]
\centering
\begin{tabular}{|c|c|c|c|}
\hline
$V(\pi_{m,n})$  &  PC $\times 10^{-6}$ (1 r.v.)&  PC$\times 10^{-6}$ (4 r.v.)\\
\hline
$\pi_{0,0}$ &  $14.711$  & $25.338$ \\\hline
$\pi_{0,1}$ &  $0.7178$ & $0.8879$ \\\hline
$\pi_{1,1}$ &  $2.4845$  & $2.5455$\\\hline
$\pi_{0,2}$ &  $0.4445$ & $0.4905$ \\\hline
$\pi_{1,2}$ &  $4.2497$  & $5.0989$ \\\hline
$\pi_{0,3}$ &  $0.7925$ & $0.8793$ \\\hline
$\pi_{1,3}$ &  $0.6677$ & $0.8450$  \\\hline
$\pi_{0,4}$ &  $1.4357$  & $1.5550$\\\hline
$\pi_{1,4}$ &  $0.9757$ & $1.5432$ \\\hline
$\pi_{0,5}$ &  $1.2546$  & $1.3592$ \\\hline
$\pi_{1,5}$ &  $3.4476$  & $4.3111$ \\\hline
\end{tabular}
\vspace{0.2cm}
\caption{\small Variance of the stationary distribution in M/M/1/5 queue 
with server breakdowns and threshold-based recovery policy}
\label{vari}
\end{table}

\begin{table}[!h]
\centering
\begin{tabular}{|c|c|c|c|}
\hline
${\rm Skew}(\pi_{m,n})$ & PC (1 r.v.)&   PC (4 r.v.) \\
\hline
$\pi_{0,0}$ &  $0.0350$ & $0.0395$   \\\hline
$\pi_{0,1}$ &  $0.0443$ & $0.0415$   \\\hline
$\pi_{1,1}$ & $0.0350$ & $0.0371$ \\\hline
$\pi_{0,2}$ & -$0.0722$ & -$0.0411$ \\\hline
$\pi_{1,2}$ &  $0.0200$  & $0.0362$ \\\hline
$\pi_{0,3}$ &  -$0.0536$  & -$0.0523$ \\\hline
$\pi_{1,3}$ &  -$0.0316$ & -$0.0556$ \\\hline
$\pi_{0,4}$ &  -$0.0011$ & -$ 0.0042$ \\\hline
$\pi_{1,4}$ &  -$0.0622$ & -$0.0116$ \\\hline
$\pi_{0,5}$ &  $0.0488$ & $0.0576$\\\hline
$\pi_{1,5}$ &  $0.0374$ & $0.0668$ \\\hline
\end{tabular}
\vspace{0.2cm}
\caption{\small Skewness coefficient of the steady-state vector in M/M/1/5 queue 
with server breakdowns and threshold-based recovery policy}
\label{skew2}
\end{table}

\begin{table}[!h]
\centering
\begin{tabular}{|c|c|c|c|}
\hline
${\rm Kurt}(\pi_{m,n})$ & PC (1 r.v.)&   PC (4 r.v.) \\
\hline
$\pi_{0,0}$  &  $3.0016$ & $3.0020$ \\  \hline
$\pi_{0,1}$ &   $3.0037$ & $3.0046$ \\    \hline
$\pi_{1,1}$ &   $3.0016$ & $3.0024$ \\      \hline
$\pi_{0,2}$ &   $3.0071$ & $3.0064$   \\  \hline
$\pi_{1,2}$ &   $2.9999$ & $3.0023$   \\  \hline
$\pi_{0,3}$ &   $3.0041$  & $3.0035$    \\   \hline
$\pi_{1,3}$ &   $3.0017$ & $3.0060$  \\    \hline
$\pi_{0,4}$ &  $2.9983$ & $2.9979$  \\     \hline
$\pi_{1,4}$ &  $3.0044$ & $ 2.9988$ \\   \hline
$\pi_{0,5}$ &  $3.0013$ & $3.0028$ \\  \hline
$\pi_{1,5}$ &  $2.9995$ & $3.0052$  \\  \hline
\end{tabular}
\vspace{0.2cm}
\caption{\small Kurtosis coefficient of the steady-state vector in M/M/1/5 queue 
with server breakdowns and threshold-based recovery policy}
\label{kur2}
\end{table}

\section{Conclusion}\label{sec5}

In this paper, we have developped a numerical approach based on polynomial
chaos expansion, to study the sensitivity and the propagation of the
epistemic uncertainty in queueing models that occurs with unreliable
servers. To illustrate the applicability of the proposed approach, two
models of queueing systems have been investigated. In the first model (M/G/1/N queue with breakdowns and repairs), 
the epistemic uncertainty only affects one input
parameter, whereas in the second model (M/M/1/N queue with server breakdowns and threshold-based recovery policy),
it affects four
parameters. In the latter case, a sensitivity analysis using Sobol' indices
is performed. When considering the stationary distribution as output
quantity of interest, it was shown that the parameter $\lambda $ of the
Poisson law modelling the customers arrival at the queue, was the most
influencial factor. This finding was then confirmed by considering $\lambda $
as the only random input parameter and by setting the three remaining ones
to their average values. In that case, it was shown that the statistical
moments of the output measure were  relatively insensitive to the other
parameters. Finally, comparisons with Monte-Carlo simulations showed the
good convergence properties of the chaos expansion.

\vspace{10pt}

\end{document}